\documentclass[11pt]{amsart}
\usepackage[]{amsmath}
\usepackage{amsfonts}
\usepackage{ifthen}
\usepackage[dvipsone]{graphics}
\usepackage{multicol}
\usepackage{caption}
\usepackage{comment}
\usepackage{color}
\usepackage{multirow}
\usepackage{amsthm}
\usepackage{fancyhdr}
\usepackage{tikz}
\usepackage{amsmath,amssymb, bbm}
\usepackage{color, colortbl}
\usepackage{amsthm}
\usepackage{graphicx}

\title{Source Code for Computing Giambelli's Formula for Type $E$ Peterson Varieties}
\author{Elizabeth Drellich}

\begin{document}
\maketitle
This document is a companion to the paper ``Monk's Rule and Giambelli's Formula for Peterson Varieties of All Lie Types."  We provide the source code for computing the Giambelli's formula in types $E_6$ $E_7$ and $E_8$.  This was done in Sage 5.2.  \\
\\
First we define a function {\tt root\_finder} that will return the root $\mathbf{r}(\mathbf{i},w)$ for a fixed reduced word $w$.\\
\\
{\tt
def T(i,z):\\
\indent    return alpha[int(z.reduced\_word()[int(i-1)])]\\
def p(i,z):\\
\indent    list=[]\\
\indent    for j in Sequence([0..int(i-2)]):\\
\indent\indent        list.append(z.reduced\_word()[int(j)])\\
\indent    return W.from\_reduced\_word(list)\\
def root\_finder(m,n):\\
\indent    a=m.reduced\_word();\\
\indent    q=n\\
\indent    for i in Sequence[0..int(m.length()-1)]:\\
\indent\indent        q=q.simple\_reflection(a[int(m.length()-i-1)]);\\
\indent    return q;\\
def r(i,w):\\
\indent    if i<=w.length():\\
\indent\indent        return root\_finder(p(i,w),T(i,w))\\
\indent    else:\\
\indent\indent        return 1\\}
\\
\noindent
These functions can be used in any Lie type.  To calculate the lists $L_{w_K}$ and $L_{p(w_K)}$ we must specify which Lie type we want to work in.  The code for calculating in type $E_8$ is given first with annotations, followed by the similar code for types $E_7$ and $E_6$.

\section{Type $E_8$}
\noindent
First we need to define the objects we will need for the computation.
{\tt
n=8\\
W=WeylGroup(['E',n],prefix="s")\\
R=RootSystem(['E',n]);\\
S=R.root\_space();\\
B=S.basis();\\
space=R.root\_lattice();\\
alpha=space.simple\_roots();\\
Q=R.root\_poset()\\
\lbrack s1,s2,s3,s4,s5,s6,s7,s8\rbrack =W.simple\_reflections()\\
w=W.long\_element()\\
\\}
\noindent
Now that we have defined the root poset on the positive roots $\Phi^+$ a simple function will give the value of $\frac{1}{t}\pi(\alpha)$ for any positive root $\alpha$.\\
\\
{\tt
def height(m):\\
\indent    return Q.rank(m)+1;\\
\\}
\noindent
The lists $L_{w_K}$ and $L_{p(w_K)}$ are obtained using this code:\\
\\
{\tt
worde8=w.reduced\_word()\\
listpwk=[]\\
for i in Sequence[1..w.length()]:\\
\indent        listpwk.append(height(r(i,w)))\\
print worde8, listpwk}\\
\\
\noindent
The list {\tt worde8} is $L_{w_{E_8}}$ and the list {\tt listpwk} is $L_{p(w_{E_8})}$. For each simple reflection $s_1,s_2,\ldots,s_8$ the evaluation of $p_{s_i}(w_{E_8})$ is calculated using this code.  Here we have used $i=3$, but this must be run eight times, setting $i$ equal to one through eight.\\
\\
{\tt 
i=3\\
q=0\\
for j in Sequence(0..(len(worde8)-1)):\\
\indent if worde8[j]==i:\\
\indent \indent q=q+listpwk[j];\\
q}\\
\\
The final value is {\tt q}$=\frac{1}{t}p_{s_i}(w_{E_8})$.  The last component of Giambelli's formula is evaluating $p_{v_K}(w_K)$.  This code finds all sublists of {\tt worde8} that are reduced words for $v_K$.  Since ${v_K}$ must end in $s_8$ and it last occurs in the $57^{th}$ spot in the list, we only need sublists of the first $57$ terms. The CPU time was slightly over $10000$ minutes for this step of the calculation.  The same calculation can be done with fewer lines of code, but a longer run time.\\
\\
{\tt
X=Set([0..56])\\
Y=X.subsets(8)\\
Z=[]\\
for q in Y:\\
\indent    if worde8[sorted(q)[7]]==8:\\
\indent\indent        if worde8[sorted(q)[6]]==7:\\
\indent\indent\indent            if worde8[sorted(q)[5]]==6:\\
\indent\indent\indent\indent                if worde8[sorted(q)[4]]==5:\\
\indent\indent\indent\indent\indent                    if worde8[sorted(q)[3]]==4:\\
\indent\indent\indent\indent\indent\indent                        if worde8[sorted(q)[2]]==3:\\
\indent\indent\indent\indent\indent\indent\indent                            if worde8[sorted(q)[1]]==2:\\
\indent\indent\indent\indent\indent\indent\indent\indent                                if worde8[sorted(q)[0]]==1:\\
\indent\indent\indent\indent\indent\indent\indent\indent\indent                                    Z.append(q)\\
\\
for q in Y:\\
\indent    if worde8[sorted(q)[7]]==8:\\
\indent\indent        if worde8[sorted(q)[6]]==7:\\
\indent\indent\indent            if worde8[sorted(q)[5]]==6:\\
\indent\indent\indent\indent                if worde8[sorted(q)[4]]==5:\\
\indent\indent\indent\indent\indent                    if worde8[sorted(q)[3]]==4:\\
\indent\indent\indent\indent\indent\indent                        if worde8[sorted(q)[2]]==2:\\
\indent\indent\indent\indent\indent\indent\indent                            if worde8[sorted(q)[1]]==3:\\
\indent\indent\indent\indent\indent\indent\indent\indent                                if worde8[sorted(q)[0]]==1:\\
\indent\indent\indent\indent\indent\indent\indent\indent\indent                                    Z.append(q)\\
\\
for q in Y:\\
\indent    if worde8[sorted(q)[7]]==8:\\
\indent\indent        if worde8[sorted(q)[6]]==7:\\
\indent\indent\indent            if worde8[sorted(q)[5]]==6:\\
\indent\indent\indent\indent                if worde8[sorted(q)[4]]==5:\\
\indent\indent\indent\indent\indent                    if worde8[sorted(q)[3]]==4:\\
\indent\indent\indent\indent\indent\indent                        if worde8[sorted(q)[2]]==3:\\
\indent\indent\indent\indent\indent\indent\indent                            if worde8[sorted(q)[1]]==1:\\
\indent\indent\indent\indent\indent\indent\indent\indent                                if worde8[sorted(q)[0]]==2:\\
\indent\indent\indent\indent\indent\indent\indent\indent\indent                                    Z.append(q)\\
}\\
Now the list {\tt Z} contains all subwords of $w_K$ that are reduced words for $v_K$ and we can sum over them as follows:\\
\\
{\tt
q=0\\
for j in Sequence(0..(len(Z)-1)):\\
\indent a=1\\
\indent for k in Sequence(0..7):\\
\indent\indent a=a*k\\
\indent q=q+a\\
q}\\
\\
This value {\tt q} is equal to $\frac{1}{t^8}p_{v_K}(w_K)$.

\section{Type $E_7$}
\noindent
The basic setup:\\
{\tt
n=7\\
W=WeylGroup(['E',n],prefix="s")\\
R=RootSystem(['E',n]);\\
S=R.root\_space();\\
B=S.basis();\\
space=R.root\_lattice();\\
alpha=space.simple\_roots();\\
Q=R.root\_poset()\\
\lbrack s1,s2,s3,s4,s5,s6,s7\rbrack =W.simple\_reflections()\\
w=W.long\_element()\\
\\
\\}
{\tt
def height(m):\\
\indent    return Q.rank(m)+1;\\
\\}
\noindent
The lists $L_{w_K}$ and $L_{p(w_K)}$ are obtained using this code:\\
\\
{\tt
worde7=w.reduced\_word()\\
listpwk=[]\\
for i in Sequence[1..w.length()]:\\
\indent        listpwk.append(height(r(i,w)))\\
print worde7, listpwk}\\
\\
\noindent
The list {\tt worde7} is $L_{w_{E_7}}$ and the list {\tt listpwk} is $L_{p(w_{E_7})}$. For each simple reflection $s_1,s_2,\ldots,s_7$ the evaluation of $p_{s_i}(w_{E_7})$ is calculated using this code.  Here we have used $i=3$, but this must be run seven times, setting $i$ equal to one through seven.\\
\\
{\tt 
i=3\\
q=0\\
for j in Sequence(0..(len(worde7)-1)):\\
\indent if worde7[j]==i:\\
\indent \indent q=q+listpwk[j];\\
q}\\
\\
The final value is {\tt q}$=\frac{1}{t}p_{s_i}(w_{E_7})$.  The last component of Giambelli's formula is evaluating $p_{v_K}(w_K)$.  This code finds all sublists of {\tt worde8} that are reduced words for $v_K$.  Since ${v_K}$ must end in $s_7$ which last appears in the $27^{th}$ spot in the list, we only need sublists of the first $27$ terms.
\\
{\tt
X=Set([0..26])\\
Y=X.subsets(7)\\
Z=[]\\
for q in Y:\\
\indent        if worde7[sorted(q)[6]]==7:\\
\indent\indent            if worde7[sorted(q)[5]]==6:\\
\indent\indent\indent                if worde7[sorted(q)[4]]==5:\\
\indent\indent\indent\indent                    if worde7[sorted(q)[3]]==4:\\
\indent\indent\indent\indent\indent                        if worde7[sorted(q)[2]]==3:\\
\indent\indent\indent\indent\indent\indent                            if worde7[sorted(q)[1]]==2:\\
\indent\indent\indent\indent\indent\indent\indent                                if worde7[sorted(q)[0]]==1:\\
\indent\indent\indent\indent\indent\indent\indent\indent                                    Z.append(q)\\
\\
for q in Y:\\
\indent        if worde7[sorted(q)[6]]==7:\\
\indent\indent            if worde7[sorted(q)[5]]==6:\\
\indent\indent\indent                if worde7[sorted(q)[4]]==5:\\
\indent\indent\indent\indent                    if worde7[sorted(q)[3]]==4:\\
\indent\indent\indent\indent\indent                        if worde7[sorted(q)[2]]==2:\\
\indent\indent\indent\indent\indent\indent                            if worde7[sorted(q)[1]]==3:\\
\indent\indent\indent\indent\indent\indent\indent                                if worde7[sorted(q)[0]]==1:\\
\indent\indent\indent\indent\indent\indent\indent\indent                                    Z.append(q)\\
\\
for q in Y:\\
\indent        if worde7[sorted(q)[6]]==7:\\
\indent\indent            if worde7[sorted(q)[5]]==6:\\
\indent\indent\indent                if worde7[sorted(q)[4]]==5:\\
\indent\indent\indent\indent                    if worde7[sorted(q)[3]]==4:\\
\indent\indent\indent\indent\indent                        if worde7[sorted(q)[2]]==3:\\
\indent\indent\indent\indent\indent\indent                            if worde7[sorted(q)[1]]==1:\\
\indent\indent\indent\indent\indent\indent\indent                                if worde7[sorted(q)[0]]==2:\\
\indent\indent\indent\indent\indent\indent\indent\indent                                    Z.append(q)\\
}\\
Now the list {\tt Z} contains all subwords of $w_K$ that are reduced words for $v_K$ and we can sum over them as follows:\\
\\
{\tt
q=0\\
for j in Sequence(0..(len(Z)-1)):\\
\indent a=1\\
\indent for k in Sequence(0..6):\\
\indent\indent a=a*k\\
\indent q=q+a\\
q}\\
\\
This value {\tt q} is equal to $\frac{1}{t^7}p_{v_K}(w_K)$.

\section{Type $E_6$}
\noindent Basic setup:\\
{\tt
n=6\\
W=WeylGroup(['E',n],prefix="s")\\
R=RootSystem(['E',n]);\\
S=R.root\_space();\\
B=S.basis();\\
space=R.root\_lattice();\\
alpha=space.simple\_roots();\\
Q=R.root\_poset()\\
\lbrack s1,s2,s3,s4,s5,s6\rbrack =W.simple\_reflections()\\
w=W.long\_element()\\
\\
\\}
{\tt
def height(m):\\
\indent    return Q.rank(m)+1;\\
\\}
\noindent
The lists $L_{w_K}$ and $L_{p(w_K)}$ are obtained using this code:\\
\\
{\tt
worde6=w.reduced\_word()\\
listpwk=[]\\
for i in Sequence[1..w.length()]:\\
\indent        listpwk.append(height(r(i,w)))\\
print worde6, listpwk}\\
\\
\noindent
The list {\tt worde6} is $L_{w_{E_6}}$ and the list {\tt listpwk} is $L_{p(w_{E_6})}$. For each simple reflection $s_1,s_2,\ldots,s_6$ the evaluation of $p_{s_i}(w_{E_6})$ is calculated using this code.  Here we have used $i=3$, but this must be run six times, setting $i$ equal to one through six.\\
\\
{\tt 
i=3\\
q=0\\
for j in Sequence(0..(len(worde6)-1)):\\
\indent if worde6[j]==i:\\
\indent \indent q=q+listpwk[j];\\
q}\\
\\
The final value is {\tt q}$=\frac{1}{t}p_{s_i}(w_{E_6})$.  The last component of Giambelli's formula is evaluating $p_{v_K}(w_K)$.  This code finds all sublists of {\tt worde8} that are reduced words for $v_K$.  Since ${v_K}$ must end in $s_7$ and the last occurrence is in the $16^{th}$ spot in the list, we only need sublists of the first $16$ terms.
\\
{\tt
X=Set([0..15])\\
Y=X.subsets(6)\\
Z=[]\\
for q in Y:\\
\indent            if worde6[sorted(q)[5]]==6:\\
\indent\indent                if worde6[sorted(q)[4]]==5:\\
\indent\indent\indent                    if worde6[sorted(q)[3]]==4:\\
\indent\indent\indent\indent                        if worde6[sorted(q)[2]]==3:\\
\indent\indent\indent\indent\indent                            if worde6[sorted(q)[1]]==2:\\
\indent\indent\indent\indent\indent\indent                                if worde6[sorted(q)[0]]==1:\\
\indent\indent\indent\indent\indent\indent\indent                                    Z.append(q)\\
\\
for q in Y:\\
\indent            if worde6[sorted(q)[5]]==6:\\
\indent\indent                if worde6[sorted(q)[4]]==5:\\
\indent\indent\indent                    if worde6[sorted(q)[3]]==4:\\
\indent\indent\indent\indent                        if worde6[sorted(q)[2]]==2:\\
\indent\indent\indent\indent\indent                            if worde6[sorted(q)[1]]==3:\\
\indent\indent\indent\indent\indent\indent                                if worde6[sorted(q)[0]]==1:\\
\indent\indent\indent\indent\indent\indent\indent                                    Z.append(q)\\
\\
for q in Y:\\
\indent            if worde6[sorted(q)[5]]==6:\\
\indent\indent                if worde6[sorted(q)[4]]==5:\\
\indent\indent\indent                    if worde6[sorted(q)[3]]==4:\\
\indent\indent\indent\indent                        if worde6[sorted(q)[2]]==3:\\
\indent\indent\indent\indent\indent                            if worde6[sorted(q)[1]]==1:\\
\indent\indent\indent\indent\indent\indent                                if worde6[sorted(q)[0]]==2:\\
\indent\indent\indent\indent\indent\indent\indent                                    Z.append(q)\\
}\\
Now the list {\tt Z} contains all subwords of $w_K$ that are reduced words for $v_K$ and we can sum over them as follows:\\
\\
{\tt
q=0\\
for j in Sequence(0..(len(Z)-1)):\\
\indent a=1\\
\indent for k in Sequence(0..5):\\
\indent\indent a=a*k\\
\indent q=q+a\\
q}\\
\\
This value {\tt q} is equal to $\frac{1}{t^6}p_{v_K}(w_K)$.
\end{document}